\documentclass{commatDV}

\usepackage[all]{xy}
\usepackage{mathtools}

\newcommand{\C}{{\mathbb C}}
\newcommand{\R}{{\mathbb R}}
\newcommand{\ZZ}{{\mathbb Z}}
\newcommand{\GG}{{\mathbb G}}

\newcommand{\cM}{{\mathcal M}}
\newcommand{\cN}{{\mathcal{N}}}
\newcommand{\cZ}{{\mathcal{Z}}}

\newcommand{\hg}{{\mathfrak{h}}}
\newcommand{\kg}{{\mathfrak{k}}}
\newcommand{\pg}{{\mathfrak{p}}}

\DeclareMathOperator{\Ad}{Ad}
\DeclareMathOperator{\Aut}{Aut}
\DeclareMathOperator{\diag}{diag}
\DeclareMathOperator{\GL}{GL}
\DeclareMathOperator{\Lie}{Lie}
\DeclareMathOperator{\Nm}{Nm}
\DeclareMathOperator{\Rr}{R}
\DeclareMathOperator{\SL}{SL}
\DeclareMathOperator{\SO}{SO}
\DeclareMathOperator{\Sp}{Sp}

\newcommand{\Dbar}{{\overline D}}
\newcommand{\GmC}{{\GG_{{\rm m},\C}}}
\newcommand{\GmR}{{\GG_{{\rm m},\R}}}
\newcommand{\ii}{{\boldsymbol{i}}}
\newcommand{\into}{\hookrightarrow}
\newcommand{\isom}{{\stackrel{\sim}{\longrightarrow}}}
\newcommand{\SMatrix}[1]{%
    \text{{\small\arraycolsep=0.4\arraycolsep\ensuremath%
    {\begin{pmatrix}#1\end{pmatrix}}%
    }}%
    }

\title[Galois cohomology of reductive algebraic groups
over the field of real numbers]%
{Galois cohomology of reductive algebraic groups\\
over the field of real numbers}

\author{Mikhail Borovoi}

\affiliation{Tel Aviv University, Israel
\email{borovoi@tauex.tau.ac.il}%
}



\abstract{
We describe functorially the first Galois cohomology set  $H^1(\R,G)$ of a connected reductive
algebraic group $G$ over the field $\R$ of real numbers in terms of a certain action of the Weyl group
on the real points of order dividing 2 of the maximal torus containing a maximal compact torus.

This result was announced with a sketch of proof in the author's 1988 note \cite{Bo}.
Here  we give a detailed proof and a few examples.
}

\VOLUME{30}
\NUMBER{3}
\YEAR{2022}
\firstpage{191}
\DOI{https://doi.org/10.46298/cm.9298}

\begin{paper}

{\em To the memory of Arkady L\!'vovich Onishchik}

\section{Introduction}
Let $G$ be a connected reductive algebraic group over the field $\R$
of real numbers.
  We wish to compute the first Galois cohomology set
$H^1(\R,G)=H^1(\text{Gal}(\C/\R), G(\C))$.
In terms of Galois cohomology one can state answers to many natural questions;
see Serre  \cite{S}, Section III.1, and  Berhuy \cite{Be}.

The Galois cohomology of classical groups and adjoint groups
is well known. The Galois cohomology of compact groups was computed
by Borel and Serre \cite{BS}, Theorem 6.8;
see also Serre's book \cite{S}, Section III.4.5.
  Here  we consider the case of a general connected reductive group over $\R$.
  We describe $H^1(\R,G)$ in terms of a certain action of the Weyl group
on the first Galois cohomology of the maximal torus containing a maximal compact torus.
Our main result  is Theorem \ref{thm:bijectivity}.

Our description of $H^1(\R,G)$ is  inspired by Borel and Serre \cite{BS}.
Our result was announced in \cite{Bo}; here we give a detailed proof and a few examples.

Since it was announced in \cite{Bo}, our Theorem \ref{thm:bijectivity} has been used in a few articles,
in particular, in \cite{Sch},  \cite{C}, and  \cite{NP}.
In  \cite{BGR}, Gornitskii, Rosengarten, and the author described, using Theorem \ref{thm:bijectivity},
the Galois cohomology of {\em quasi-connected} reductive $\R$-groups
(normal subgroups of connected reductive $\R$-groups).
Our description in  \cite{BGR} is similar to that of Theorem \ref{thm:bijectivity}.
In \cite{BE} Evenor and the author used Theorem \ref{thm:bijectivity}
to describe {\em explicitly} the Galois cohomology of simply connected semisimple $\R$-groups.
In \cite{BT} and \cite{BT-bis}, Timashev and the author used Theorem \ref{thm:bijectivity}
to describe explicitly the Galois cohomology of connected reductive $\R$-groups.

Note that cited articles refer to Theorem 9 of an early preprint version of this note.
In this published version, Theorem 9 became Theorem \ref{thm:bijectivity}.

\section{Preliminaries}
We recall the definition of the first Galois cohomology set $H^1(\R,G)$
of an algebraic group $G$ defined over $\R$.
  The set of 1-cocycles is defined by
$Z^1(\R,G)=\{z\in G(\C)\ |\ z\bar{z}=1\}$
where the bar denotes  complex conjugation.
  The group $G(\C)$ acts on the right on $Z^1(\R,G)$ by
\[
z*x=x^{-1}z\bar{x},
\]
where $z\in Z^1(\R,G)$ and $x\in G(\C)$.
By definition $H^1(\R,G)= Z^1(\R,G)/G(\C)$.
Let $G(\R)_2$ denote the subset of elements of $G(\R)$ of order 2 or 1.
Then  $G(\R)_2\subset Z^1(\R,G)$, and we obtain a canonical map
$G(\R)_2\to H^1(\R,G)$.

\begin{lemma}\label{lem2.1}
Let $S$ be an algebraic $\R$-torus. Let $S_0$ denote the largest
compact (that is, anisotropic) $\R$-subtorus in $S$, and let $S_1$ denote
the largest split subtorus in $S$. Then :
\begin{enumerate}
\item[\rm (a)] The map $\lambda\colon S(\R)_2\to H^1(\R,S)$ induces a canonical
isomorphism
\[
    S(\R)_2/S_1(\R)_2 \isom H^1(\R,S).
\]

\item[\rm(b)] The composite map $\mu\colon S_0(\R)_2\to H^1(\R,S_0)\to H^1(\R,S)$ is surjective.

\item[\rm(c)]$(S_0\cap S_1)(\R)=S_0(\R)_2\cap S_1(\R)_2$, and the surjective map
$\mu$   of (b) induces an isomorphism
$$S_0(\R)_2/(S_0\cap S_1)(\R)\isom H^1(\R,S).$$
\end{enumerate}
\end{lemma}

\begin{proof}
Any $\R$-torus is isomorphic to a direct product of tori of three types,
see Casselman \cite{Cass}, Section 2:

(1) $\GmR$,

(2) $\Rr_{\C/\R}\GmC$,

(3) $\Rr^1_{\C/\R}\GmC$.
  Here  $\GG_{\rm m}$ denotes the multiplicative group, $\Rr_{\C/\R}$ denotes the
Weil restriction of scalars, and
\begin{equation*}
\Rr^1_{\C/\R}\GmC=\ker\big[{\rm Nm}_{\C/\R}\colon \Rr_{\C/\R}\GmC\to\GmR\big],
\end{equation*}
where ${\rm Nm}_{\C/\R}$ is the norm map.

We prove (a).
The composite homomorphism $S_1(\R)_2\into S(\R)_2\to H^1(\R,S)$ factors via $H^1(\R,S_1)=1$, and hence it is trivial.
We obtain an induced homomorphism
 \[S(\R)_2/S_1(\R)_2 \to H^1(\R,S);\]
we must prove that it is bijective.
It suffices to consider the three cases:

(1) $S=\GmR$, that is, $S(\R)=\R^\times$.
Then  $H^1(\R,S)=1$.
We have $S_1=S$, so $S(\R)_2/S_1(\R)_2=1$. This proves (a) in case (1).

(2) $S=\Rr_{\C/\R}\GmC$, that is $S(\R)=\C^\times$.
Then  $H^1(\R,S)=1$.
We have $S_1=\GmR$, $S_1(\R)=\R^\times$,
$S_1(\R)_2=\{1,-1\}=S(\R)_2$, so $S(\R)_2/S_1(\R)_2=1$. This proves (a)
in case (2).

(3) $S=\Rr^1_{\C/\R}\GmC$, that is $S(\R)=\{x\in\C^\times\ |\ \Nm(x)=1\}$,
where $\Nm(x)=x\bar{x}$.
  Then  by the definition of Galois cohomology
$H^1(\R,S)=\R^\times/\Nm(\C^\times)\simeq\{-1, 1\}$.
The homomorphism $S(\R)_2=\{-1,1\}\to H^1(\R,S)$ is an isomorphism.
This proves (a) in case (3).

Assertion (b) reduces to the cases (1), (2), (3), where it is obvious
(note that only in case (3)  we have $H^1(\R,S)\neq 1$).

Concerning (c), we have a commutative diagram
$$
 \xymatrix{
S_0(\R)_2 \ar@{^{(}->}[r] \ar[d] \ar[rd]^{\mu}  &S(\R)_2 \ar[d]^{\lambda} \\
H^1(\R, S_0) \ar[r]           & H^1(\R,S).
}
$$
We see from (a) that $\ker \mu=S_0(\R)_2\cap S_1(\R)_2$,
and we know from (b) that $\mu$ is surjective.
Thus we obtain a canonical isomorphism
$$
S_0(\R)_2/(S_0(\R)_2\cap S_1(\R)_2)\isom H^1(\R,S).
$$
It remains only to check that $S_0(\R)_2\cap S_1(\R)_2=(S_0\cap S_1)(\R)$.
This can be easily checked in each of the cases (1), (2), (3)
(note that only in case (2) this group is nontrivial).
\end{proof}

\begin{corollary}\label{cor2.2}
Assume that $S$ is an $\R$-torus such that $S=S'\times S''$, where
$S'$ is a compact torus and $S''=\Rr_{\C/\R} T$, where $T$ is a $\C$-torus.
Then  $H^1(\R,S)=H^1(\R,S')=S'(\R)_2$.
\end{corollary}

\begin{proof}
The assertion follows from the proof of Lemma \ref{lem2.1}(a),
because $S'$ is a direct product of tori of type (3); hence
$H^1(\R,S')=S'(\R)_2$\;, and $S''$ is a direct product of tori of type (2),
whence $H^1(\R,S'')=1$.
\end{proof}

We say that a connected real algebraic group $H$ is \emph{compact},
if the group $H(\R)$ is compact, that is, $H$ is reductive and
anisotropic.
We shall need the following two standard facts.

\begin{lemma}[well-known]
\label{lem:compact}
Any nontrivial semisimple algebraic group $H$ over $\R$
contains a nontrivial connected compact subgroup.
\end{lemma}

\begin{proof}
This assertion follows from the classification,
(see, for instance, Helgason \cite{He}, Section X.6.2, Table V).
We prove it without using the classification.

Let $\kappa\colon\hg\times \hg\to \R$ denote the Killing form on $\hg$.
Let $\hg=\kg+\pg$ be a Cartan decomposition
of the real semisimple Lie algebra $\hg=\Lie\;H$.
This means that the linear transformation
\[\theta\colon \hg\to\hg,\ \  k+p\,\mapsto\, k-p\quad\text{for}\ \,k\in\kg,\, p\in\pg\]
is an automorphism of $\hg$, and that the bilinear form
\[b_\theta(x,y)=-\kappa\big(x,\theta(y)\big)\]
is positive definite on $\hg$.
Set
\[K=\big\{h\in H\mid \Ad\, h\in O(\hg,b_\theta)\big\}.\]
Then  $K$ is a real algebraic subgroup of $H$.
We have $\Lie\, K=\kg$; see Gorbatsevich, Onishchik, and Vinberg \cite{OV2}, Section 4.3.2.
Since $H(\R)$ has finitely many connected components and the center of $H(\R)^0$ is finite,
by \cite{OV2}, Corollary 5 of Theorem 4.3.2, the group $K(\R)$ is compact.

Since $\pg$ and $\kg$ are the eigenspaces of $\theta$ with eigenvalues $-1$, and $+1$, respectively,
we have  $[\pg,\pg]\subset\kg$.
  If $\kg=0$, then $\hg=\pg$ and $[\pg,\pg]=0$, whence $\hg$ is commutative,
which is clearly impossible.
Thus $\kg\neq 0$.
But $\kg$ is the Lie algebra of the identity component $K^0$ of $K$,
which is a connected compact algebraic subgroup of $H$.
Thus $H$ contains a nontrivial connected compact algebraic subgroup.
\end{proof}

\begin{lemma}[well-known] \label{lem:max-compact}
Any two maximal compact tori in a connected reductive real algebraic
group $H$ are conjugate under $H(\R)$.
\end{lemma}

\begin{proof}
It suffices to prove that any two maximal compact tori in the derived
group $[H,H]$ of $H$ are conjugate.
  This follows from the following well-known
facts from the theory of  Lie groups:
(1) Any two maximal compact subgroups in a connected semisimple Lie
group are conjugate (see, for instance, Gorbatsevich, Onishchik, and Vinberg \cite{OV2}, Section 4.3.4, Theorem 3.5);
(2) Any two maximal tori in a connected compact Lie group are conjugate
(see, for instance, Onishchik and Vinberg \cite{OV1}, Section 5.2.7, Theorem 15).
\end{proof}

\section{Main result}
Let $G$ be a connected reductive algebraic group over $\R$.
  Let $T_0$ be a maximal compact torus in $G$.
  Set $T=\cZ(T_0)$, $N_0=\cN(T_0)$, $W_0=N_0/T$, where $\cZ$ and $\cN$
denote the centralizer and the normalizer in $G$, respectively.

We prove that $T$ is a torus.
  By Humphreys \cite{Hu}, Theorem 22.3 and Corollary 26.2.A, the centralizer
$T$  of $T_0$ is a connected reductive $\R$-group.
  The torus $T_0$ is a maximal compact torus in $T$, and it is central
in $T$.
   Since by Lemma \ref{lem:max-compact}
all the maximal compact tori in $T$ are conjugate under $T(\R)$,
we see that $T_0$ is the only maximal compact torus in $T$.
  It follows that the derived group $[T,T]$ of $T$
contains no nontrivial compact tori.
  By Lemma \ref{lem:compact} every nontrivial semisimple group over $\R$ has a nontrivial
compact connected algebraic subgroup, hence a nontrivial compact torus.
  We conclude that $[T,T]=1$, and hence $T$ is a torus.
  We see that $T$ is a {\em fundamental torus in $G$}, that is,
a maximal torus containing a maximal compact torus.

We have a right action of $W_0$ on $T_0$ defined by
$(t,w)\mapsto t\cdot w\coloneqq n^{-1}tn$, where $t\in T_0(\C)$, $n\in N_0(\C)$,
$n$ represents $w\in W_0(\C)$. This action is defined over $\R$.
  We prove that $W_0(\C)$ acts on $T_0$ effectively.
  Indeed, if $w\in W_0(\C)$ with representative $n\in N_0(\C)$ acts
trivially on $T_0$, then $n^{-1}tn=t$ for any $t\in T_0(\C)$, and
hence $n\in T(\C)$ (because the centralizer of $T_0$ is $T$), whence $w=1$.

We prove that $W_0(\C)=W_0(\R)$.
 We have seen that $W_0(\C)$ embeds in $\Aut_\C(T_0)$.
 Since $T_0$ is a compact torus, all the complex automorphisms of $T_0$
are defined over $\R$.
  We see that the complex conjugation  acts trivially on $\Aut_\C(T_0)$,
and hence on $W_0(\C)$.
  Thus $W_0(\R)=W_0(\C)$.

Note that $N_0$ normalizes $T$; hence $W_0$ acts on $T$.
We define a right action $*$ of $W_0(\R)$ (which is equal to $W_0(\C)$) on $H^1(\R,T)$.
  Let $z\in Z^1(\R,T)$, $n\in N_0(\C)$, $z$ represents $\xi\in H^1(\R,T)$,
$n$ represents $w\in W_0(\R)=W_0(\C)$.
  We set
\begin{equation*}
\xi * w=[n^{-1}z\bar{n}]=[n^{-1}zn\cdot n^{-1}\bar{n}],
\end{equation*}
where brackets  $[\ \ ]$ denote the cohomology class.

We prove that $*$ is a well defined action.
  First, since $N_0$ normalizes $T$ and $z\in T(\C)$, we see
that $n^{-1}zn\in T(\C)$.
  Now $w\in W_0(\R)$, whence $w^{-1}\bar{w}=1$ and
$n^{-1}\bar{n}\in T(\C)$.
  It follows that $n^{-1}z\bar{n}=n^{-1}zn\cdot n^{-1}\bar{n}\in T(\C)$.
We have
\begin{equation*}
n^{-1}z\bar{n}\cdot \overline{n^{-1}z\bar{n}}=
      n^{-1}z\bar{n}\bar{n}^{-1}\bar{z}n=1
\end{equation*}
because $z\bar{z}=1$.
  Thus $n^{-1}z\bar{n}\in Z^1(\R, T)$.
   If $z'\in Z^1(\R,T)$ is another representative of $\xi$,
then $z'=t^{-1}z\bar{t}$ for some $t\in T(\C)$, and
\begin{equation*}
n^{-1}z'\bar{n}=n^{-1}t^{-1}z\bar{t}\bar{n}
      =(n^{-1}tn)^{-1}\cdot n^{-1}z\bar{n}\cdot \overline{n^{-1}tn}
      =(t')^{-1}(n^{-1}z\bar{n})\overline{t'}
\end{equation*}
where $t'=n^{-1}tn$, $t'\in T(\C)$.
   We see that the cocycle $n^{-1}z'\bar{n}\in Z^1(\R,T)$
is cohomologous to $n^{-1}z\bar{n}$.
   If $n'$ is another representative of $w$ in $N_0(\C)$, then
$n'=nt$ for some $t\in T(\C)$, and
  $(n')^{-1}z\overline{n'}=t^{-1}n^{-1}x\bar{n}\bar{t}$.
  We see that $(n')^{-1}z\overline{n'}$ is cohomologous to
$n^{-1}z\bar{n}$.
Thus $*$ is indeed a well defined action
of the group $W_0(\R)$ on the set $H^1(\R,T)$.

  Note that in general $[1]*w=[n^{-1}\bar{n}]\neq [1]$, and therefore,
the action $*$  does not respect the group structure in $H^1(\R,T)$.

Let $\xi\in H^1(\R,T)$ and $w\in W_0(\R)$.
  It follows from the definition of the action $*$  that the images of
$\xi$ and $\xi * w$ in $H^1(\R,G)$ are equal.
  We see that the map $H^1(\R,T)\to H^1(\R,G)$ induces a map
$H^1(\R,T)/W_0(\R)\to H^1(\R,G)$.

The following theorem is the main result of this note:

\begin{theorem}\label{thm:bijectivity}
  Let $G$, $T_0$, $T$, and $W_0$ be as above.
The map
$$H^1(\R,T)/W_0(\R)\to H^1(\R,G)$$
induced by the map $H^1(\R,T)\to H^1(\R,G)$ is a bijection.
\end{theorem}

\begin{proof} We prove the surjectivity.
  It suffices to show that the map $H^1(\R,T)\to H^1(\R,G)$ is surjective.
This was proved by Kottwitz \cite{Ko}, Lemma 10.2,
with a reference to Shelstad \cite{Sh}.
  We give a different proof.
  Let $\eta\in H^1(\R,G)$, $\eta=[z]$, $z\in G(\C)$, $z\bar{z}=1$.
  Let $z=us=su$, where $s$ and $u$ are the semisimple and the unipotent
parts of $z$, respectively (see Humphreys \cite{Hu}, Theorem 15.3).
  We have $us\bar{u}\bar{s}=1$, where $\bar{u}\bar{s}=\bar{s}\bar{u}$
(because $us=su$).
 Thus $us=\bar{u}^{-1}\bar{s}^{-1}$, where $u$ and $\bar{u}^{-1}$
are unipotent, $s$ and $\bar{s}^{-1}$ are semisimple, $us=su$.
From the equality $\bar{u}\bar{s}=\bar{s}\bar{u}$ it follows that
$\bar{u}^{-1}\bar{s}^{-1}=\bar{s}^{-1}\bar{u}^{-1}$.
  Since the Jordan decomposition in $G(\C)$ is unique
(see Humphreys \cite{Hu}, Theorem 15.3),
we conclude that $s=\bar{s}^{-1}$, $u=\bar{u}^{-1}$.
  In other words, $s\bar{s}=1$, $u\bar{u}=1$, that is, $s$ and $u$ are cocycles.

Since $u$ is unipotent, the logarithm $\log(u)\in\Lie\; G_\C$ is defined.
We have:
\[
\log(u)+\overline{\log(u)}=0.
\]
Set $y=\frac{1}{2} \log(u)$, then $y+\bar{y}=0$.
  We have $-y+\log(u)+\bar{y}=0$, where $-y$, $\bar{y}$ and $\log(u)$
pairwise commute.
  Set $u'=\exp(y)$, then $(u')^{-1} u \overline{u'}=1$.
Since $s$ commutes with $u$, we have $\Ad(s)y=y$, and hence $s$ commutes with
$u'$.
  We obtain $(u')^{-1}su\overline{u'}=s$, and hence the cocycle $z=su$ is cohomologous
to the cocycle $s$, where $s$ is semisimple.

We may and shall therefore assume that $z$ is semisimple.
  Set $C=\cZ_{G_\C}(z)$.
  Since $\bar{z}=z^{-1}$, we have $\overline{C}=C$, and hence the algebraic subgroup
$C$ of $G_\C$ is defined over $\R$.
  The semisimple element $z$ is contained in a maximal torus of $G_\C$
(see Humphreys \cite{Hu}, Theorem 22.2); hence $z$ is contained in
the identity component $C^0$ of $C$.
  The group $C^0$ is reductive, see Steinberg \cite{St}, Section 2.7(a).
  Let $T'$ be a maximal torus of $C^0$ defined over $\R$, then
$z\in T'(\C)$, because $z$ is contained in the center of $C^0$.
  By Lemma \ref{lem2.1}(b) the class $\eta$ of $z$ comes from
the maximal compact subtorus $T'_0$ of $T'$.
  By Lemma \ref{lem:max-compact}
any compact torus in $G$ is conjugate under $G(\R)$ to a subtorus of $T_0$.
Thus $\eta$ comes from $H^1(\R,T_0)$, hence from $H^1(\R,T)$.
 This proves the surjectivity in Theorem \ref{thm:bijectivity}.

We prove the injectivity in  Theorem \ref{thm:bijectivity}.
  Let $z,z'\in T(\C)$, $z\bar{z}=1$, $z'\overline{z'}=1$,
$z=x^{-1}z'\bar{x}$, where $x\in G(\C)$.
  We shall prove that $z=n^{-1}z'\bar{n}$ for some $n\in N_0(\C)$.

For $g\in G(\C)$ set $g^\nu=z\bar{g}z^{-1}$.
  Then  $\nu$ is an involutive antilinear automorphism of $G_\C$,
and in this way we obtain a twisted form ${}_z G$ of $G$.
  Since $z\in T(\C)$, the embeddings of the tori $T_\C$ and $T_{0,\C}$
into ${}_zG_\C$ are defined over $\R$.
We denote the corresponding $\R$-tori of $_z G$ again by $T$ and $T_0$, respectively.
  The centralizer of $T_0$ in ${}_zG$ is $T$.
 The compact torus $T_0$ of $_z G$ is contained in some maximal compact torus $S$ of $_z G$,
and clearly $S$ is contained in the centralizer $T$ of $T_0$ in $_z G$.
Since $T_0$ is the largest compact subtorus of $T$, we conclude that the $S=T_0$.
Thus $T_0$ is a maximal compact torus in ${}_z G$.

Consider the embedding
$i_x\colon t\mapsto x^{-1} t x\colon T_{0,\C}\to {}_z G_\C$\,.
  We have $i_x(t)^\nu=z\bar{x}^{-1}\bar{t}\bar{x}z^{-1}$.
  Since $z\bar{x}^{-1}=x^{-1}z'$, we obtain
\begin{equation*}
z\bar{x}^{-1}\bar{t}\bar{x}z^{-1}=x^{-1}z'\bar{t}(z')^{-1}x=x^{-1}\bar{t}x
    =i_x(\bar{t})\,.
\end{equation*}
We see that $i_x(t)^\nu=i_x(\bar{t})$; hence $i_x$ is defined over $\R$.
  Set $T_0'=i_x(T_0)$; it is a compact algebraic torus in ${}_z G$,
and $\dim T_0'=\dim T_0$.
  Therefore, the torus $T_0'$ is conjugate to $T_0$ under ${}_z G(\R)$,
say, $T_{0,\C}=h^{-1}T_{0,\C}' h$, where $h\in{}_z G(\R)$.
Set $n=xh$.
  Then
\[n^{-1}T_{0,\C}\,n=h^{-1}x^{-1}T_{0,\C}\,xh=h^{-1}T_{0,\C}'h=T_{0,\C},\]
whence $n\in N_0(\C)$.
  The condition $h\in{}_z G(\R)$ means that $z\bar{h}z^{-1}=h$, or
$h^{-1}z\bar{h}=z$.
  It follows that
\begin{equation*}
n^{-1}z'\bar{n}=h^{-1}x^{-1}z'\bar{x}\bar{h}=h^{-1}z\bar{h}=z.
\end{equation*}
We have proved that there exists $n\in N_0(\C)$ such that
$z=n^{-1}z'\bar{n}$, and hence the cohomology classes
$[z],[z']\in H^1(\R,T)$ lie in the same orbit of $W_0(\R)$ in $H^1(\R, T)$.
This proves the injectivity in Theorem \ref{thm:bijectivity}.
\end{proof}

\begin{remark}
If $G$ is a compact group, then Theorem \ref{thm:bijectivity}
asserts that
\[
    H^1(\R,G)=T(\R)_2/W,
\]
where $T$ is a maximal torus in
$G$, and $W$ is the Weyl group with the usual action.
  This was earlier proved by Borel and Serre \cite{BS}.
\end{remark}

\begin{remark}
The real form $G$ of $G_\C$ defines an involutive automorphism $\tau$
of the {\em based root datum of $G_\C$}, see \cite{BT}, Proposition 3.7,
and hence an involutive automorphism $\tau_D$ of the Dynkin diagram $D=D(G_\C)$.
This automorphism $\tau_D$ is trivial if and only if the derived group $[G,G]$ of $G$
is an {\em inner form} of a compact group, that is, has a compact maximal torus.
Let $\Dbar$ denote the {\em twisted Dynkin diagram} corresponding to $D$ and $\tau_D$.
Then  $W_0$ is isomorphic to the Weyl group of $\Dbar$; see \cite{BT}, Proposition 7.11(iii).
This Coxeter group is described in the book of Carter \cite{Carter}, Chapter 13.
\end{remark}

\section{Examples}

In this section, written following a suggestion of the referee, we compute, {\em using Theorem~\ref{thm:bijectivity}},
the sets $H^1(\R,G)$ when $G=$ $\GL_{n,\R}$, $\Sp_{2m,\R}$, $\SO_{p,q}$.

\begin{example}
Let $G=\GL_{2m,\R}$, \,$m\in\ZZ_{>0}$.
For $z=a+b\ii\in \C$, we write
\[ \cM(z)=\SMatrix{a&b\\-b&a}.\]
 Consider the tori $T$ and  $T_0$ in $G$ such that
\begin{align*}
&T(\R)=\big\{\diag(\cM(z_1),\dots,\cM(z_m))\ \,\big|\  \, z_k=a_k+b_k\ii\in\C^\times,\ k=1,\dots, m \big\},\\
&T_0(\R)=\big\{\diag(\cM(z_1),\dots,\cM(z_m))\in T(\R)\ \,\big|\ \,z_k=a_k+b_k\ii,\ a_k^2+b_k^2=1\big\}.
\end{align*}
Then  $T_0$ is a maximal compact torus in $G$, and $T$ is a fundamental torus containing $T_0$.
Since $T\simeq (\Rr_{\C/\R}\GmC)^m$, by the proof of Lemma \ref{lem2.1}, case (2), we have $H^1(\R,T)=\{1\}$,
and by Theorem \ref{thm:bijectivity} we conclude that $H^1(\R,G)=\{1\}$.

Similarly, if $G=\GL_{2m+1,\R}$, then $G$ has a fundamental torus
\[T\simeq\GmR\times (\Rr_{\C/\R}\GmC)^m.\]
Again we have  $H^1(\R,T)=\{1\}$ and $H^1(\R,G)=\{1\}$.
Note that it is well known that $H^1(K,\GL_n)=\{1\}$ for any $n$ and any field $K$;
see Serre \cite{S-LF}, Section X.1, Proposition 3.
\end{example}

\begin{example}
Let $G=\SL_{2,\R}$.
It has a maximal torus $T=T_0$ with group of $\R$-points
\[T(\R)=\big\{\cM(z)\ \big|\ z=a+b\ii,\ a^2+b^2=1\big\}. \]
Set $n=\diag(\ii,-\ii)\in G(\C)$; then $n\in N(\C)$, where $N=N_0=\cN_G(T)$.
We have $\#H^1(\R,T)=2$ with representatives $\cM(1),\cM(-1)$. An easy calculation shows that
\[ n^{-1} \cM(-1)\bar n= n^{-1} \cM(-1) n\cdot n^{-1}\bar n=\cM(1)=1.\]
Thus $H^1(\R,T)/W_0=\{1\}$, and by Theorem \ref{thm:bijectivity} we have $H^1(\R,G)=1$.

Note that in this case the group $W=W_0=N_0/T$ has two elements,
and that the element $w\coloneqq[n^{-1}]\in W(\R)$ has no representative in $N_0(\R)$,
because otherwise we would have $[1]*w=[1]$.
Thus in this case $W_0(\R)\neq N_0(\R)/T(\R)$.
\end{example}

\begin{example}
Let $G=\Sp_{2m}=\Sp(\R^{2m},\psi)$, where $\psi$ is the skew-symmetric bilinear form with matrix
\[M_\psi=\diag(J, \dots,J)\quad\text{where}\ J=\cM(\ii)=\SMatrix{0&1\\-1 &0}.\]
The group $G$ has a compact maximal torus $T=T_0$ with
\[T(\R)=\big\{\diag(\cM(z_1),\dots,\cM(z_m))\ \,\big|\ \,z_k=a_k+b_k\ii,\ a_k^2+b_k^2=1\big\}.\]
We have
\[T(\R)_2=\big\{\diag(\cM(z_1),\dots,\cM(z_m))\in T(\R)\ \,\big|\ \,z_k=\pm 1\big\}.\]
Let $t=\diag(\cM(z_1),\dots,\cM(z_m))\in T(\R)_2$.
Write
\[n=\diag(n_1,\dots,n_m),\quad\text{where}\ n_k=
\begin{cases}
\diag(\ii,-\ii) &\text{if}\ z_k=-1,\\
\diag(1,1)      &\text{if}\ z_k=1.
\end{cases}
\]
Then  $n\in N_0(\C)$ and
\[ n^{-1}\cdot t\cdot \bar n=1.\]
We see that for any $[t]\in H^1(\R,T)$ there exists $w=[n]\in W_0(\R)$ with
$$[t]*w=[1].$$
Thus $H^1(\R,T)/W_0=\{1\}$, and by Theorem \ref{thm:bijectivity} we have $H^1(\R,G)=\{1\}$.
Note that it is well known that $H^1(K,\Sp_{2m})=\{1\}$ for any $m$ and any field $K$;
see Serre \cite[Section~III.1.2, Proposition~3]{S} .
\end{example}

\begin{example}
Let $G=\SO(p',p'')=\SO(\R^{p'+p''}, f)$, where $f$ is the diagonal quadratic form with matrix
\[ M_{ f}=\diag\big(\underbrace{+1,\dots,+1}_{p'\text{ times}},\underbrace{-1,\dots,-1}_{p''\text{ times}}\big).\]
We consider the case when both $p'$ and $p''$ are even: $p'=2r',\ p''=2r''$.
Our group $G$ has a compact maximal torus $T=T_0$ with group of $\R$-points
\[T(\R)=\big\{\diag\big(\cM(z_1),\dots\cM(z_{r'+r''})\big)\ \,\big|\ \, z_k=a_k+\ii b_k,\ a_k^2+b_k^2=1\big\}.\]
We have
\[T(\R)_2=\big\{\diag\big(\cM(z_1),\dots\cM(z_{r'+r''})\big)\ \,\big|\ \, z_k=\pm1\big\}.\]
The Weyl group $W=W(G_\C,T_\C)$ is isomorphic to $(\pm1)^{r'+r''-1}\rtimes S_{r'+r''}$,
where $S_{r'+r''}$ is the symmetric group on the $r'+r''$ symbols $1,\dots,r'+r''$.
The subgroup $(\pm1)^{r'+r''-1}$ acts on $T(\R)_2$ trivially.
We compute $T(\R)_2/S_{r'+r''}$.

For a subset $\Xi\subseteq\{1,\dots,r'+r''\}$, we set
\[c_\Xi=\diag\big(\cM(z_1),\dots,\cM(z_{r'+r''})\big)\ \text{with}\ z_k=
\begin{cases}
-1 &\text{if}\ k\in\Xi,\\
+1 &\text{otherwise}.
\end{cases}
\]
Then
\[T(\R)_2=\big\{c_\Xi\ \big|\ \Xi\subseteq \{1,\dots,r'+r''\}\big\}.\]

Consider the subgroup $S_{r'}\times S_{r''}\subseteq S_{r'+r''}$,
where $S_{r''}$ is the symmetric group on the $r''$ symbols $r'+1,\dots,r'+r''$.
Its elements are represented by elements of $N(\R)$,
 and hence they act on $T(\R)_2$ by the usual conjugation:
\[c_\Xi*\sigma=c_{\sigma^{-1}\Xi} \quad\text{for}\ \,\sigma\in S_{r'}\times S_{r''}.\]
Write $\Xi=\Xi'\cup\Xi''$ where
\[\Xi'=\Xi\cap \{1,\dots,r'\},\quad \Xi''=\Xi\cap\{r'+1,\dots,r'+r''\}.\]
We see that the $W$-orbit of $c_\Xi$ depends only on the cardinalities of $\Xi'$ and $\Xi''$.

The group $S_{r'+r''}$ is generated by its subgroup $S_{r'}\times S_{r''}$
and $\sigma_{1,r'+1}=(1,r'+1)$.
In order to compute the action of $\sigma_{1,r'+1}$ on $T(\R)_2$,
we consider the case $r'=1$, $r''=1$, $G=\SO(2,2)$.
Consider the block matrix
\[n=\ii\SMatrix{0&I_2\\ I_2&0}\quad\text{where}\ \, I_2=\diag(1,1).\]
One can check that $n\in N(\C)\subset G(\C)$ and $n$ represents $\sigma_{1,2}$.
Let
\[c=c_{\{1,1\}}=\diag(-1,-1,-1,-1).\]
We have
\[n^{-1}c\bar n=cn^{-1}\bar n=1.\]

Returning to the case of arbitrary $r'$ and $r''$,
we see that when both $\Xi'$ and $\Xi''$ are non-empty,
an element of $\Xi'$ can be cancelled with an element of $\Xi''$.
Thus the $W$-orbit of a cocycle $c_\Xi$ depends
only on the difference  $\#\Xi'-\#\Xi''$.
For $s'$ and $s''$ such that  $1\le s'\le r'$,\, $1\le s''\le r''$,  we write
\begin{align*}
&c'_{s'}=c_{\Xi'}\quad\text{with}\ \ \Xi'=\{1,\dots,s'\},\\
&c''_{s''}=c_{\Xi''}\quad \text{with}\ \ \Xi''=\{r'+1,\dots,r'+s''\}.
\end{align*}
By Theorem \ref{thm:bijectivity} we conclude that
\[\#H^1(\R,G)=r'+r''+1\]
with representatives
\[1\,\cup\,\big\{c'_{s'}\mid 1\le s'\le r'\big\}\,\cup\,\big\{c''_{s''}\mid 1\le s''\le r''\big\}.\]

The cases $G=\SO(2r',2r''+1)$ and $G=\SO(2r'+1,2r''+1)$ are similar to the case $\SO(2r',2r'')$;
in both cases we have $\#H^1(\R,G)=r'+r''+1$.

Alternatively, one can use the fact that $H^1(\R, G)$ for $G=\SO(\R^n,f)$
classifies isomorphism classes of real quadratic forms $f'$ on $\R^n$
with $\det M_{ f'}=\det M_{ f}$ (see Serre \cite{S-LF}, Section X.2, Proposition 4),
and one can classify the isomorphism classes of such $f'$ using Sylvester's law of inertia.

\end{example}

\subsection*{Acknowledgements}
The author was partially supported by the Hermann Minkowski Center for Geometry.
The author thanks the anonymous referee for helpful suggestions,
and  Dmitry A. Timashev for helpful email correspondence.

\end{paper}